\documentclass{article}
\usepackage{graphicx,url,amsthm}

\usepackage{geometry}
\usepackage{amssymb,latexsym,enumitem,bbm,amsmath}
\usepackage[noadjust]{cite}
\usepackage{setspace}
\usepackage{color}
\usepackage{hyperref}
\usepackage{bbm}

\newcommand{\R}{\mathbb{R}}
\newcommand{\N}{\mathbb{N}}

\newcommand{\1}{\mathbbm{1}}
\newcommand{\assign}{:=}

\newcommand{\ito}{It\^o }

\numberwithin{equation}{section}

\newtheorem{theorem}{Theorem}[section]

\newtheorem{lemma}[theorem]{Lemma}
\newtheorem{proposition}[theorem]{Proposition}
\newtheorem{remark}[theorem]{Remark}

\newtheorem{definition}[theorem]{Definition}

\newtheorem{example}[theorem]{Example}

\providecommand{\keywords}[1]
{
  \small	
  \textbf{\textit{Keywords---}} #1
}

\providecommand{\subjclass}[1]
{
  \small	
  \textbf{\textit{Mathematics Subject Classification 2020---}} #1
}

\allowdisplaybreaks
\title{Self-normalized scaled quadratic variation}

\author{
	\textsc{Suprio Bhar}
	\thanks{Department of Mathematics and Statistics, Indian Institute of Technology Kanpur, India (E-mail: \it{suprio@iitk.ac.in})}
	\and
    \textsc{Purba Das}
	\thanks{Department of Mathematics, King's College London, UK (E-mail: \it{purba.das@kcl.ac.uk})}
	\and
	\textsc{Barun Sarkar} 
	\thanks{Department of Mathematics, Indian Institute of Technology Madras, India (E-mail: \it{barun@iitm.ac.in})}
}
\date{}

\begin{document}

\maketitle
\begin{abstract}
\noindent 
The concept of a scaled quadratic variation was originally introduced by E. Gladyshev in 1961 for processes with Gaussian increments. Using certain deterministic scaling, arrived at from the covariance of the process, Gladyshev showed that the sum of scaled square increments along the dyadic partition sequence converges almost surely to a finite limit. In this paper, we propose a pathwise counterpart in which the deterministic normalization is replaced by a self-normalizing factor built from the $p$-th variation of the path along a given sequence of partitions. The resulting quantity requires no probabilistic assumption and no knowledge of a covariance structure, and its scale is both path-dependent and sensitive to the
partition sequence. Under a mild regularity condition on the limiting $p$-th variation, we show that the self-normalized and the classical deterministic normalizations are comparable, and for fractional Brownian motion the two agree up to a multiplicative constant. We establish a switching behaviour in the index, and prove that for $p \ge 2$ the self-normalized scaled quadratic variation obeys a smooth-transformation formula under $C^2$ maps; at $p=2$ this recovers the known transformation rule for quadratic variation. Since only squared increments are scaled, the construction polarizes, yielding a matrix-valued scaled quadratic variation for every $p \ge 1$ for $\R^d$ valued paths.
We conclude with examples beyond the Gaussian setting.

\end{abstract}


    
\section{Introduction}
In 1956, Baxter \cite{Baxter1956} showed that for any Gaussian process $X_t$ for which the variance of increment is of order $\bigtriangleup t$, the sum of squared increments over a refining partition converges almost surely to the integrated variance. Gladyshev \cite{gladyshev1961} then generalized this result of Baxter for Gaussian processes $X_t$ beyond the variance structure. Gladyshev only assumed that the Gaussian process $X_t$ needs to follow  $\mathbb E|X(t+h)-X(t)|^2 \sim h^{2H}$ (no stationarity of increments or self-similarity was assumed). For such $X_t$, Gladyshev showed that along dyadic sequence of partitions, a scaled version of the square increment (and hence the name `scaled quadratic variation')
\begin{align}\label{eq: scaledQV}
    \sum_{i=1}^{2^n} 2^{-n(\gamma-1)} \left(X\left(\frac{i+1}{2^n}\right)- X\left(\frac{i}{2^n}\right)\right)^2
\end{align}
converges almost surely to an `increment variance' type term for a very particular choice of $\gamma$. Because \eqref{eq: scaledQV} involves only squared increments, it is attractive computationally and extends to vector-valued paths without difficulty.
\par The quantity $\gamma$ is closely linked to the order of roughness of the underlying process $X_t$, it comes from the scaling factor of the underlying Gaussian process $X_t$, and is linked with the covariance. For example, for a self-similar Gaussian process with Hurst exponent $H$, there is exactly one exponent $\gamma = 2-2H$ for which \eqref{eq: scaledQV} has a finite non-zero limit almost surely and for any $\gamma \neq 2 -2H$ the limit is either $0$ or $+\infty$.
\par The asymptotic behaviour of scaled quadratic variations of stochastic processes has been studied extensively since Gladyshev's seminal paper in 1961, well beyond the Gaussian setting (for example, semimartingales). In particular for semimartingales, the almost-sure limit theory of normalized power variations has been extensively explored by Jacod and Protter \cite{jacod-protter-discretization}, and for more general self-similar
processes similar type of limit result has been studied in the monograph of Tudor \cite[Part II]{tudor2013}, including a treatment of the non-Gaussian case. This type of scaled quadratic variation was later used for estimation of the Hurst exponent in various contexts (see, for example, \cite{leahy-EJP25, kubilius2010, han-Schied2025,mishura2017-book}). What we would like to stress here is that, likewise to the original scaled quadratic variation introduced in \eqref{eq: scaledQV}, the normalization in all of the follow-up results is \emph{deterministic}: the scale is a power of the mesh $\bigtriangleup t$, or of $2^{-n}$. However, if one wants to develop the limit theorem of scaled increments for a deterministic path $X:[0, T]\to \R$, and not necessarily for a process with known covariance structure like before, one needs to come up with a normalizing/scaling factor which does not depend/ is inspired from the covariance structure of the process. This is precisely what we explore in this manuscript. 

\par To this end, we first summarise how roughness is currently captured from a purely pathwise perspective for any continuous path $X$. The concept of limits of square increment from a pure pathwise perspective (no probability involved) was introduced 1981 by F\"ollmer \cite{follmer1981}, extending the definition beyond semimartingales/ self-similar processes. For any function/ process $X$ with `similar roughness' to semimartingales, F\"ollmer showed that the (pathwise) quadratic variation of $X$ along a large class of partition sequences (denoted as $[X]_\pi$ for a partition sequence $\pi$) is finite (this can be thought of pathwise analogue of Baxter \cite{Baxter1956}). Nonetheless, processes like fractional Brownian motion (fBM) with Hurst exponent $H\neq \frac{1}{2}$, and differential equations driven by fBM, exhibit strictly rougher or strictly smoother behaviour than semimartingales. For these processes/functions, the quadratic variation is either $0$ or $\infty$ and hence (pathwise) quadratic variation fails to measure their variability. To solve this issue, in the spirit of F\"ollmer \cite{follmer1981}, for all $1<p<\infty$, the concept of $p$-th variation along a sequence of partitions $\pi$, denoted by $[X]^{(p)}_\pi$ (and a corresponding It\^o type formula)  was introduced by Cont \& Perkowski in \cite{perkowski2019} (see Definition \ref{def:p-var} below). It turns out, likewise to scaled quadratic variation, the notion of $p$-th variation along a given partition sequence also satisfies a switching behaviour. i.e. there exists a critical variation index $p^*$ (potentially dependent on the choice of partition $\pi$) such that for all $p\neq p^*$, the $p$-th variation along a sequence of partitions $\pi$ is either $0$ or $\infty$. However, as Definition \ref{def:p-var} shows, the $p$-th variation is defined as the sum of absolute increments of $X$ raised to power $p$, there is no natural extension to vector-valued paths for $p\not\in 2\N$. To this end, we develop a pathwise version of scaled quadratic variation of \eqref{eq: scaledQV}. One of the main contributions of this manuscript is to identify the correct (purely pathwise) scaling factor in \eqref{eq: scaledQV}. The idea we used in Definition \ref{def:scaled QV gen} is that, instead of using the covariance structure of the path $x$ to define the weights, we use the $p$-th variation of the path itself as a normalizing factor.

\par Our self-normalized scaled quadratic variation is bounded above and below by constant multiples of the classical scaled quadratic variation under some nominal additional (path) regularity condition on $x$. Furthermore, both quadratic variation and $p$-th variation along a partition sequence $\pi$ satisfy a `smooth-transformation formula' (see \cite{ananova2017,perkowski2019,ruhong2024}) governed by $|f'|^p$ as follows: 
\begin{align*}
   [f\circ X]^{(p)}_\pi(t) = \int_0^t |f'(X(u))|^p \; d [ X]^{(p)}_\pi(u) \qquad \text{ for any `smooth' function } f.
\end{align*}
It turns out our path-adaptive/ self-normalized scaled quadratic variation (see Theorem \ref{Thm.ito.isomerty} for details) also satisfies a similar `smooth-transformation formula'. Finally, it turns out that the self-normalized scaled quadratic variation also satisfies a similar switching behaviour, i.e. there exists a critical variation index $p^*$ (potentially dependent on the choice of partition $\pi$) such that for all $p\neq p^*$, the $p$-th variation along a sequence of partitions $\pi$ is either $0$ or $\infty$ (see Lemma \ref{lem: switching} below).

{\it Preview:} This paper is organised as follows. In section \ref{sec:2} we introduce the notations, relevant definitions and basic properties involving self-normalized scaled quadratic variation. Section \ref{sec:3} introduces a pathwise smooth transformation-type formula for (self-normalized) scaled quadratic variation in the spirit of F\"ollmer, and Section \ref{sec:4} concludes the paper with some examples with finite non-trivial scaled quadratic variation, including a two-dimensional fractional Brownian motion.

\section{Definitions and properties}\label{sec:2}
\subsection{Definitions}
Here we introduce, in the spirit of F\"ollmer \cite{follmer1981}, the concept of $p$-th variation along a (deterministic) sequence of partitions $\pi^n=(t_0^n, \dots, t^n_{N(\pi^n)})$ with $0=t_0^n<...< t^n_k<...< t^n_{N(\pi^n)}=T$ (See e.g. \cite{perkowski2019} and the references therein for details). We denote by $N(\pi^n)$, the number of intervals in the partition $\pi^n$ and
\begin{align*}
     |\pi^n|=\sup\{ |t^n_i-t^n_{i-1}|, i=1,\cdots,N(\pi^n)\},\qquad\underline{\pi^n }=\inf\{ |t^n_i-t^n_{i-1}|, i=1, \cdots , N(\pi^n)\},
\end{align*} 
the size of the largest (resp. the smallest) interval of $\pi^n$. Throughout the paper, unless mentioned otherwise, we take $x \in C([0,T],\R)$, and write
\begin{align}\label{eq:del}
    \bigtriangleup t^n_i := t^n_{i+1}-t^n_i\;, \qquad \bigtriangleup x(t^n_i):= x(t^n_{i+1}) - x(t^n_i).
\end{align}

\noindent Define the \emph{oscillation} of $x \in C([0,T],\R)$ along $\pi_n$ as
\[
	\mathbf{osc}(x,\pi_n) := \max_{[t_j,t_{j+1}] \in \pi_n} \max_{r,s \in [t_j,t_{j+1}]} |x(s) - x(r)|.
\]
Here and in the following we write $[t_j,t_{j+1}] \in \pi_n$ to indicate that $t_j$ and $t_{j+1}$ are both in $\pi_n$ and are immediate successors (i.e. $t_j < t_{j+1}$ and $ \pi_n \cap(t_j , t_{j+1}) = \emptyset$).

\begin{definition}[$p$-th variation along a sequence of partitions]\label{def:p-var} Let $p\geq 1$. A  continuous path $x \in C([0,T],\R)$ is said to have a $p$-th variation along a sequence of
partitions $\pi=(\pi_n)_{n\geq 1}$ if $\mathbf{osc}(x,\pi_n)\to 0$ and the sequence of measures
\begin{align}\label{Eq: p-th var}
    \mu^n \assign \sum_{[t_j, t_{j + 1}] \in \pi_n} \delta (\cdot - t_j) |\bigtriangleup x (t_j)|^p
\end{align}	
converges weakly to a measure $\mu$ without atoms. In that case we write $x \in V^p_\pi([0,T],\R)$ and $[x]_\pi^{(p)}(t) := \mu([0,t])$ for $t \in [0,T]$, and we call $[x]_\pi^{(p)}$ the \emph{$p$-th variation} of $x$. We shall also write $[x]_{\pi_n}^{(p)}(t)$ for $\mu^n([0, t])$ for $t \in [0, T]$.
\end{definition}

As a consequence of $\mu$ having no atoms one get $x\in V^p_\pi([0,T],\R)$ implies $t\mapsto [x]_\pi^{(p)}(t) $ is continuous. The same conclusion also holds for Definition \ref{def:scaled QV gen} below. 
\par Note that this notion of $p$-th variation is different from $p$-variation in the sense of sup over partitions.
The Gladyshev estimator \cite{gladyshev1961} was originally introduced in the context of Gaussian processes, and there have been many follow-up works in this setup \cite{leahy-EJP25,MR2965922}. We now propose a (pathwise) version of the scaled quadratic variation that does not require any Gaussian assumptions about the underlying process. The convergences are also understood in a uniform sense instead of a probabilistic $L^2$ sense. 
\par Throughout the manuscript we will use the convention $0^0 = 1$. Furthermore, we adopt the convention $0^\gamma \cdot \1_{>0} := 0$ 
for all $\gamma \in \mathbb{R}$, so that intervals on which the $p$-th variation increment vanishes contribute nothing to the final sum, regardless of the sign of $\gamma = \frac{p-2}{p}$. More precisely, for any $a \geq 0$ 
and $\gamma \in \mathbb{R}$, we define:
\begin{align}\label{eq: convention}
    a^{\gamma} \cdot \1_{a > 0} := 
    \begin{cases} 
        a^{\gamma} & \text{if } a > 0 \\ 
        0          & \text{if } a = 0
    \end{cases}.
\end{align}

\begin{remark}
For $x\in V^p_\pi([0,T],\R)$ and $[t_i^n,t_{i+1}^n]\in \pi_n$, we shall write 
\[\Delta[x]_\pi^{(p)} (t^n_i) := [x]_\pi^{(p)}(t_{i+1}^n)-[x]_\pi^{(p)}(t_i^n)\]
or, more simply, $\Delta[x]_\pi^{(p)} (t_i)$ by suppressing $n$ if the corresponding partition $\pi_n$ is clear from the context. Note that the increments of a continuous path $x$ or the variable $t$ are being denoted by $\bigtriangleup$, and the increments of any $p$-th variation by $\Delta$ for distinction.
\end{remark}

Let us now define the notion of (self-normalized) scaled quadratic variation.
\begin{definition}\label{def:scaled QV gen} Let $p\geq1$ and $\gamma=\frac{p-2}{p}$. A  continuous path $x \in V^p_\pi([0,T],\R)$ is said to have scaled quadratic variation\footnote{Throughout the paper, whenever the index $p$ is clear from the context, we just say that the path has finite scaled quadratic variation.} of index $p$ along a sequence of
partitions $\pi=(\pi_n)_{n\geq 1}$, if $\mathbf{osc}(x,\pi_n)\to 0$ and the sequence of measures
\begin{align}\label{eq: mu^n} 
    \mu^n_{SQV} (\cdot) \assign \sum_{[t_j, t_{j + 1}] \in \pi_n} \delta (\cdot - t_j) \left(\Delta[x]_\pi^{(p)} (t_j)\right)^\gamma|\bigtriangleup x (t_j)|^2  \1_{\Delta[x]_\pi^{(p)} (t_j)>0},  
\end{align}
converges weakly to a measure $\mu_{SQV}$ without atoms. In that case we write $x \in \mathbb{S}^p_\pi ([0,T],\R)$ and $\langle x\rangle^{(p)}_\pi(t) := \mu_{SQV}([0,t])$ for $t \in [0,T]$, and we call $\langle x\rangle^{(p)}_\pi(t)$ the \emph{scaled quadratic variation} of $x$ at time $t$. We will sometimes denote the finite quantity $\mu^n_{SQV}$ as $ \langle x\rangle^{(p)}_{\pi^n}$.
\end{definition}

Note that in Definition \ref{def:scaled QV gen} the same $p$ is used while taking the increment of $p$-th variation and in the exponent $\gamma$. This is not a random choice; rather, this is the only interesting choice. Let us define a more general definition with $p,q > 0$ and $x \in V^p_\pi([0,T],\mathbb{R})$. We define the 
scaled quadratic variation of order $q$ with reference $p$ along $\pi$ as the 
limit (when it exists) of the measures
\[
    \mu^{p,q}_{n,\pi}(\cdot) := \sum_{[t_j,t_{j+1}]\in\pi_n} \delta(\cdot - t_j)
    \Bigl(\Delta[x]_\pi^{(p)} (t_j)\Bigr)^{\frac{q-2}{q}}
    (\bigtriangleup x (t_j))^2 \mathbbm{1}_{\Delta[x]_\pi^{(p)} (t_j)>0}.
\]
Then for all $t\in[0,T]$ the limiting measure $\mu^{p,q}_{\pi}([0,t])$ satisfies switching behaviour, i.e for $q<p$; $\mu^{p,q}_{\pi}([0,t])=\infty $ and for $q>p$; $\mu^{p,q}_{\pi}([0,t])=0$. We explain this in detail in Lemma \ref{lem: switching}.

For a higher-dimensional path $x = (x^1, \cdots, x^d)\in C([0,T],\R^d)$, the quadratic variation of $x$ along a partition sequence $\pi$ is defined as a $[x]_\pi \in C([0,T ], S^+_d)$, where $S^+_d$ is the set of $d\times d$ symmetric positive semidefinite matrices.
For the case of $p$-th variation with $p\in 2\N$, the p-th variation can be defined in $\mathbb S_p(\R^d)$, which is defined as direct sum of symmetric $k$-tensor $Sym_k(\R^d)$ for $k =0, 1, 2,\cdots,p$ (See \cite[Section 4]{perkowski2019}). However, for $p\not\in 2\N$, such a higher-dimensional analogue is not present. i.e. the function $t\mapsto [x]^{(p)}_\pi(t)$ is not defined \cite[Definition 4.1.]{ruhong2024}. 

One crucial advantage of our scaled quadratic variation (Definition \ref{def:scaled QV gen}) is that, since we only rely on the square of the increment of the underlying path $x$, we can introduce a $d\times d$ matrix-valued scaled quadratic variation for all $p\geq 1$, following the approach of F\"ollmer for quadratic variation. In particular, this allows us to define a scaled quadratic covariation term for all $p\geq 1$.  

\par Recall the usual quadratic variation of a path $x$ at level $\pi^n$ is given by the sum of $(\Delta x(t^n_i))^2$. Similarly, the usual quadratic covariation of two paths $x,y$ at level $\pi^n$ is given by the sum of $(\Delta x(t^n_i) )\times (\Delta y(t^n_i))$, which allows the quadratic variation to be defined for higher-dimensional paths. 
We will introduce a notion of scaled quadratic variation and scaled quadratic covariation by simply replacing the increments of $x$ and $y$ by a renormalized/ scaled increment at each level $\pi^n$. To this end, for one dimensional path $x,y\in V^p_\pi([0,T],\R )$, let us define the renormalized/ scaled increment at level $\pi^n$ as   
\[ \delta_i^p x := \left(\Delta [x]_\pi^{(p)}(t_i) \right)^{\frac{p-2}{2p}} (\bigtriangleup x (t_i))\, \1_{\Delta[x]_\pi^{(p)} (t_i)>0} \qquad \delta_i^p y := \left(\Delta [y]_\pi^{(p)}(t_i) \right)^{\frac{p-2}{2p}} (\bigtriangleup y (t_i)) \, \1_{\Delta[y]_\pi^{(p)} (t_i)>0} 
\]
We will now define the scaled quadratic covariation of $x$ and $y$ at level $\pi^n$ as sum of $\delta_i^p x \times  \delta_i^p y$.
\begin{definition}[Scaled quadratic covariation]\label{def:scaled_covariation} Consider a sequence of
partitions $\pi=(\pi_n)_{n\geq 1}$. Let $p\geq 1$ and let  $x=(x^1,\ldots,x^d)\in V_\pi^p([0,T],\mathbb R^d)$. Set $\gamma:=\frac{p-2}{p}$. For $k,\ell\in\{1,\ldots,d\}$ define 
\begin{align*}
\mu^{n,(k,\ell)}_{SQV} := \sum_{[t_j,t_{j+1}]\in\pi_n} \delta(\cdot-t_j) \Big( \Delta[x^k]_\pi^{(p)} (t_j) \,\Delta[x^\ell]_\pi^{(p)} (t_j) \Big)^{\gamma/2} (\bigtriangleup x^k (t_i))\,&(\bigtriangleup x^\ell (t_j)) \\ &\times \mathbbm{1}_{\Delta[x^k]_\pi^{(p)} (t_j) \,\Delta[x^\ell]_\pi^{(p)}(t_j) >0} .\end{align*}
If $\mu^{n,(k,\ell)}_{SQV}$ converges weakly to a measure $\mu_{\pi}^{(k,\ell)}$ without atoms, then we define \[ \langle x^k,x^\ell\rangle_\pi^{(p)}(t) := \mu_{\pi}^{(k,\ell)}([0,t]), \qquad t\in[0,T],\]    
and call it the scaled quadratic covariation of order $p$ of $x^k$ and $x^\ell$ along $\pi$. The matrix-valued process 
\[ t\mapsto \langle x\rangle_\pi^{(p)} (t):= \big( \langle x^k,x^\ell\rangle_\pi^{(p)} (t)\big)_{k,\ell=1}^d \] is called the scaled quadratic variation of the $d$ dimensional path $x$. \end{definition}
\begin{remark}
The broader goal of introducing (pathwise) scaled quadratic variation is to understand the path properties of the underlying process for varying roughness index $p$ and identify the correct scaling factor to normalise the $p$-th variation. Since the roughness/variation index is closely related to the finiteness of $p$-th variation (\cite{das2022a}), we want to get inspiration for the self-normalizing factor from the $p$-th variation itself. The following calculations show that the $ p$-th variation is indeed closely linked to the proposed pathwise scaled quadratic variation with the same index $p$. Hence, we choose $\mathbb S^p_\pi([0,T],\R) \subset V^p_\pi([0,T],\R)$ with the same (roughness) index $p$. Consider $x$ to be a continuous function with finite $p$-th variation along $\pi$. Then for all $t\in[0,T]$ the $p$-th variation along $\pi$ can be rewritten as:
\begin{align*}
[x]^{(p)}_{\pi}(t) & = \lim_{n\to\infty} \sum_{\pi^n\cap [0,t]} |\bigtriangleup x (t^n_i)|^p  
\\& = \lim_{n\to\infty} \sum_{\pi^n\cap [0,t]} |\bigtriangleup x (t^n_i)|^{p-2} \1_{|\bigtriangleup x (t^n_i)| > 0}\times |\bigtriangleup x (t^n_i)|^2
\\& = \lim_{n\to\infty} \sum_{\pi^n\cap [0,t]} |[x]_{\pi^n}^{(p)}(t^n_{i+1})-[x]_{\pi^n}^{(p)}(t^n_{i})|^\gamma \1_{\Delta [x]_{\pi^n}^{(p)}(t^n_{i}) > 0}\times |\bigtriangleup x (t^n_i)|^2 \qquad \text{where, } \gamma = \frac{p-2}{p}.
\end{align*}
The scaled quadratic variation (Definition \ref{def:scaled QV gen}) is obtained simply by replacing $\pi^n$ with $\pi$ in the first increment. Note that for $p<2$ the exponent $\gamma <0$. So in case $\Delta[x]^{(p)}_\pi(t^n_i) = 0$ in some interval $[t^n_i,t^n_{i+1}]$, the normalizing factor is undefined, hence the convention in \eqref{eq: convention} is introduced in the definition. 
\par Furthermore, in the spirit of Freedman's construction \cite{freedman,davis2018}, for any continuous $x$ one can construct a partition sequence $\tilde\pi$ with vanishing mesh such that the $p$-th variation along this partition sequence is equal to $0$ (See \cite[Lemma~2.7 \& Proposition~2.8]{daskim2025} for the construction of $\tilde\pi$). Since in Definition \ref{def:scaled QV gen} the scaling comes from the increment of the $p$-th variation, along the same choice of partition sequence $\tilde\pi$, we also have $\langle x\rangle^{(p)}_{\tilde\pi}\equiv 0$, suggesting our definition coincides with the notion of $p$-th variation, and captures the dependence of roughness along the choice of partition sequence. It is not clear whether if one takes the $2^{-n}$ or $ \bigtriangleup t^n_{i}$ as the scaling factor in the definition of scaled quadratic variation (for example \eqref{eq: scaledQV}, as suggested in \cite{leahy-EJP25}), this would still capture the correct path properties along a given sequence of partitions.
\end{remark}

We provide a sufficient condition under which the scaled quadratic variation of index $p$ exactly matches with the $p$-th variation.
\begin{lemma}\label{SQV-recovers-pth-var}
Let $p>2$ and $\pi=(\pi_n)$ be a sequence of partitions of $[0,T]$ with vanishing mesh, and let $x\in V_\pi^p([0,T],\mathbb{R})$. Suppose that for every interval $I_i^n=[t_i^n,t_{i+1}^n]\in \pi_n$, we have 
\[ \Delta[x]_\pi^{(p)} (t^n_i) = |\bigtriangleup x (t^n_i)|^p(1+\eta_i^n), \] 
where 
\[ \max_{i\in\{1,\cdots N(\pi_n)\}} |\eta_i^n|\xrightarrow[]{n\to\infty} 0. \] 
Then the measures 
\[ \mu^n_{SQV} := \sum_{I_i^n\in\pi_n} \delta(\cdot -t_i^n) \big(\Delta[x]_\pi^{(p)}(t^n_i)\big)^\gamma (\bigtriangleup x (t^n_i))^2  \1_{\Delta[x]_\pi^{(p)} (t^n_i)>0}\] converge weakly to the Stieltjes measure $d[x]_\pi^{(p)}$. In particular, 
\[ \langle x\rangle_\pi^{(p)}(t) = [x]_\pi^{(p)}(t), \qquad t\in[0,T]. \] 
\end{lemma} 

\begin{proof} Since we assume $p>2$, we will drop the $\1_{\cdot}$ in the Definition \ref{def:scaled QV gen}, as $\gamma>0$. For notational simplicity set $u_i^n:=|\bigtriangleup x (t^n_i)|^p$. So by assumption, $\Delta[x]_\pi^{(p)}(t^n_i) = u_i^n(1+\eta_i^n)$. Hence 
\begin{align*} \big(\Delta[x]_\pi^{(p)} (t_i^n)\big)^\gamma (\bigtriangleup x (t^n_i))^2 &= \big(u_i^n(1+\eta_i^n)\big)^\gamma (u_i^n)^{2/p}  = (u_i^{n}) (1+\eta_i^n)^\gamma . \end{align*} 
Following the notation of \eqref{Eq: p-th var} and \eqref{eq: mu^n} we get,
\[ \mu^n_{SQV}-\mu^n = \sum_{i} u_i^n\Big((1+\eta_i^n)^\gamma-1\Big) \delta(\cdot -t_i^n). \] 
Therefore, 
\begin{align}\label{Eq: total var}
   \|\mu^n_{SQV}-\mu^n \|_{TV} \leq \sup_i\big|(1+\eta_i^n)^\gamma-1\big| \sum_{i} u_i^n \xrightarrow[]{n\to\infty} 0, 
\end{align}
as a consequence of  $\sup_i|\eta_i^n|\to0$ and $\sum_i u_i^n = \sum_{t^n_i<T} |\bigtriangleup x (t^n_i)|^p \xrightarrow[]{n\to\infty} [x]_\pi^{(p)}(T)$. 
Moreover, for every $t\in[0,T]$, 
\[ \mu^n([0,t]) = \sum_{\{i:t_i^n<t\}} |\bigtriangleup x (t^n_i)|^p \xrightarrow[]{n\to\infty} [x]_\pi^{(p)}(t). \]
Now $\mu^n\Rightarrow d[x]^{(p)}_\pi$ is exactly the hypothesis  $x\in V^p_\pi([0,T],\mathbb R)$. Since convergence in total variation  implies weak convergence, \eqref{Eq: total var} yields $\mu^n_{SQV}-\mu^n\Rightarrow 0$, and hence $\mu^n_{SQV}\Rightarrow d[x]^{(p)}_\pi$. The limit has no atoms, so $x\in\mathbb S^p_\pi([0,T],\mathbb R)$ with $\langle x\rangle_\pi^{(p)}(t)=[x]_\pi^{(p)}(t)$ for all $t\in[0,T]$.
\end{proof}

\begin{remark}
A. Schied \cite[Proposition~2.7]{schied2016} provides an example to show that for dyadic partition, the class $V^p_\pi([0,T], \R)$ is not a vector space even for $p=2$. In a similar line of argument, it can be shown that the class $\mathbb S^p_\pi([0,T],\R)$ is also not a vector space. 
\end{remark}

\begin{remark}\label{counter.example}
The limit in \eqref{eq: mu^n} may not exist in general. So we propose a relaxed version of scaled quadratic variation where the limit in \eqref{eq: mu^n} is replaced by either $\limsup$ or $\liminf$. The benefit of considering such a relaxed version is that the scaled quadratic variation will always exist in such cases. At the same time, one cannot have a change of variable/ stability under smooth transformation under such relaxation (see Theorem \ref{Thm.ito.isomerty}). But if one considers such a relaxation by replacing $\lim$ with $\limsup$ or $\liminf$ then it may turn out that $\infty=\limsup\langle x\rangle^{(p)}_{\pi^n} \neq\liminf \langle x\rangle^{(p)}_{\pi^n} = 0$. These types of examples can also be constructed in the $p$-th variation context.
\end{remark}
\par \noindent The following example gives a path $x\in C^0([0,1],\R)$ for which $$\infty=\limsup\langle x\rangle^{(2)}_{\pi^n}(t) =\limsup\; [x]_{\pi^n}^{(2)}(t) \neq \liminf\; [x]_{\pi^n}^{(2)}(t)=\liminf \langle x\rangle^{(2)}_{\pi^n}(t) = 0.$$
\begin{example}\label{counter.example.1}
Consider a function with its dyadic Schauder expansion as follows (please see \cite[Section 3.2]{das2021} for further details), where $\mathbb T$ represents the dyadic sequence of partitions $\mathbb T_n$ on $[0, 1]$. 
\begin{align*}
    x(t) = \lim_{n\to\infty} \sum_{m=0}^{n-1}\sum_{k=0}^{2^m-1} \theta_{m,k} e^{\mathbb T}_{m,k}(t),
\end{align*}
where,
$$\theta_{m,k} = \begin{cases} \sqrt{2\ell-\frac{\ell-1}{2^{\ell-1}}} & \text{for }m=S_\ell-1 \text{ and } \forall k \\
0 & \text{otherwise}    
\end{cases},$$
where, $S_\ell = 1+2+\cdots + \ell$ for all $\ell\in \mathbb N$.
This is a continuous function as a result of \cite[
Lemma 2.10]{bayraktar2025}.
The quadratic variation (and by convention the scaled quadratic variation) at level $n$ can be represented as (see \cite{schied2016, das2021} for details),
\begin{align*}
    [x]_{\mathbb T^n}^{(2)}(1) = \frac{1}{2^n}\sum_{m=0}^{n-1} \sum_{k=0}^{2^m-1} \theta_{m,k}^2.
\end{align*} So we get:
\begin{align*}
\langle x\rangle^{(2)}_{\mathbb T^{S_n}} (1)= [x]_{\mathbb T^{S_n}}^{(2)}(1) &= \frac{1}{2^{S_n}}\sum_{m=0}^{S_n-1} \sum_{k=0}^{2^m-1} \theta_{m,k}^2
= \frac{1}{2^{S_n}} \bigg[\sum_{i=1}^n 2^{S_i-1} \theta_{m,k}^2\bigg]
 = \cdots
 = \frac{1}{2^{S_n}} n 2^{\frac{n(n+1)}{2}} \xrightarrow[]{n\to\infty} \infty.
\end{align*}
Similarly:
\begin{align*}
\langle x\rangle^{(2)}_{\mathbb T^{S_n-1}}(1) = [x]_{\mathbb T^{S_n-1}}^{(2)}(1) = \frac{1}{2^{S_n-1}}\sum_{m=0}^{S_n-2} \sum_{k=0}^{2^m-1} \theta_{m,k}^2 &= \frac{1}{2^{S_n-1}} \bigg[\sum_{i=1}^{n-1} 2^{S_i-1} \theta_{m,k}^2\bigg] \\& = \cdots = \frac{2}{2^{S_n}} (n-1) 2^{\frac{(n-1)n}{2}} \xrightarrow[]{n\to\infty} 0.
\end{align*}
This completes the example.
\end{example}

\subsection{Properties and Lemmata}
The following lemma compares the new proposed pathwise definition of scaled quadratic variation (Definition \ref{def:scaled QV gen}) to the classical definition, defined as a limit of \eqref{eq: scaledQV}.
\begin{lemma}\label{Lemma: connection with classical theory}
Take $\pi$ to be a sequence of partitions for $[0,T]$ with vanishing mesh and $x\in V^p_\pi ([0,T],\R)$. Denote $\phi(t) = [x]^{(p)}_\pi(t)$. If $t\mapsto\phi(t)\in C^1_b([0,T],\R)$ and $\exists c,C$ such that $0<c\leq \phi'\leq C<\infty$, then for $\gamma = \frac{p-2}{p}$,
\begin{align*}
\limsup_n\langle x\rangle^{(p)}_{\pi^n}(T) <\infty &\iff \limsup_n \sum_{i=0}^{N(\pi^n)-1} |\bigtriangleup t^n_i|^\gamma |\bigtriangleup x (t^n_i)|^2 <\infty,\\
\liminf_n\langle x\rangle^{(p)}_{\pi^n}(T) >0 &\iff \liminf_n \sum_{i=0}^{N(\pi^n)-1} |\bigtriangleup t^n_i|^\gamma |\bigtriangleup x (t^n_i)|^2 >0.
\end{align*}
\end{lemma}

\begin{proof}[Proof of Lemma \ref{Lemma: connection with classical theory}]
Since, $\phi^\prime \geq c >0$, $\phi$ is strictly increasing and we drop the $\1$ from our calculation. The first equivalence can be reduced as follows.
\begin{align*}
   & \limsup_n\langle x\rangle^{(p)}_{\pi^n}(T) <\infty 
   \iff  \limsup_n \sum_{i=0}^{N(\pi^n)-1} \left(\Delta\phi(t^n_{i})\right)^\gamma|\bigtriangleup x (t^n_i)|^2  <\infty.
\end{align*}
Since $\phi$ is differentiable on compact intervals, the mean value theorem tells there exists $u^n_{i,i+1}\in [t^n_{i}, t^n_{i+1}]$ such that the above statement has the following equivalent formulation, 
\begin{align} \label{Eq: sup}
   \iff & \limsup_n \sum_{i=0}^{N(\pi^n)-1} \left(\phi'(u^n_{i,i+1})\right)^\gamma\left(\bigtriangleup t^n_i\right)^\gamma|\bigtriangleup x (t^n_i)|^2 <\infty \nonumber\\
\implies & \limsup_n \sum_{i=0}^{N(\pi^n)-1} \left(\bigtriangleup t^n_i\right)^\gamma|\bigtriangleup x (t^n_i)|^2 <\infty ,
\end{align}
where we use, $\phi'(u)^\gamma \geq \min\{c^\gamma,C^\gamma\}>0$.
From \eqref{Eq: sup} one can also get the other side of the inequality as follows. 
\begin{align*}
& \limsup_n \sum_{i=0}^{N(\pi^n)-1} \left(\bigtriangleup t^n_i\right)^\gamma|x(t^n_{i + 1}) - x(t^n_i)|^2 <\infty  \\
   \implies &  \limsup_n \sum_{i=0}^{N(\pi^n)-1}  \left( \frac{\Delta \phi(t^n_{i})}{\phi'(u^n_{i,i+1})}\right)^\gamma |\bigtriangleup x (t^n_i)|^2<\infty\\
\implies &  \limsup_n \sum_{i=0}^{N(\pi^n)-1}  \left( \frac{\Delta \phi(t^n_{i})}{C\1_{\gamma\geq 0}+ c\1_{\gamma<0}}\right)^\gamma |\bigtriangleup x (t^n_i)|^2<\infty\\
\implies &  \limsup_n \sum_{i=0}^{N(\pi^n)-1}  \left(\Delta\phi(t^n_{i})\right)^\gamma |\bigtriangleup x (t^n_i)|^2<\infty  \qquad \text {since, } 0<c<C<\infty.
\end{align*}
This completes the if and only if argument for the $\limsup_n$. The inequality involving $\liminf_n$ can be shown in a similar line of arguments. 
\end{proof}
Below, we characterise the cutoff properties of pathwise scaled quadratic variation. 

\begin{lemma}[Switching behaviour of \texorpdfstring{$\mathbb{S}^p_\pi ([0,T],\R)$}{TEXT}] \label{lem: switching}
Take $\pi$ to be a sequence of partitions for $[0,T]$ with vanishing mesh and $x\in V^p_\pi ([0,T],\R)$. Denote $\phi(t) = [x]^{(p)}_\pi(t)$. If $t\mapsto\phi(t)\in C^1_b([0,T],\R)$ and $\exists \; c,C$ such that $0<c\leq \phi'\leq C<\infty$, then the following holds. 
\begin{enumerate}
    \item[(1)] If $x\in \mathbb S^p_\pi([0,T],\R)$ with $\langle x\rangle^{(p)}_\pi(t)>0$, then
    \begin{align*}
\mu^{p,q}_{\pi}([0,t])= \begin{cases}
    \infty & \text{for } q<p\\
    0 & \text{for } q>p\\
\end{cases}.       
    \end{align*}
    
\item[(2a)]  If $\limsup_n\langle x\rangle^{(p)}_{\pi^n}(t) <\infty$, then 
$\limsup_n \mu^{p,q}_{n,\pi}([0,t])\equiv 0$ for all $q> p$.
\item[(2b)] If $\limsup_n\langle x\rangle^{(p)}_{\pi^n}(t) >0$, then $\limsup_n 
\mu^{p,q}_{n,\pi}([0,t]) \equiv \infty$ for all $q< p$.
\item[(3a)]  If $\liminf_n\langle x\rangle^{(p)}_{\pi^n}(t) <\infty$, then $\liminf_n
\mu^{p,q}_{n,\pi}([0,t])\equiv 0$ for all $q> p$.
\item[(3b)] If $\liminf_n\langle x\rangle^{(p)}_{\pi^n}(t) >0$, then $\liminf_n
\mu^{p,q}_{n,\pi}([0,t])\equiv \infty$ for all $q< p$.
\end{enumerate}
\end{lemma}

\begin{proof}
Since, $\phi^\prime \geq c >0$, $\phi$ is strictly increasing and we drop the $\1$ from our calculation.
\begin{enumerate}
 \item[(1)] Consider the case $q=p+\epsilon$, i.e. for all $q>p$, for all $t\in[0,T]$:
 \begin{align*}
 0\leq
 \mu^{p,q}_{\pi}([0,t])
 & = \lim_{n\to\infty} \sum_{[t^n_i, t^n_{i + 1}] \in \pi^n\cap[0,t]} \left( \bigtriangleup x(t^n_i)\right)^2 \left( \Delta\phi(t_i^n)\right)^{\frac{p+\epsilon-2}{p+\epsilon}}  \\
& = \lim_{n\to\infty} \sum_{\pi^n\cap[0,t]} \left( \bigtriangleup x(t^n_i)\right)^2 \phi^\prime(\eta_{i,i+1})^{\frac{q-2}{q}}\left( \bigtriangleup t_i^n\right)^{1-\frac{2}{p}} \left( \bigtriangleup t_i^n\right)^{\frac{2}{p}-\frac{2}{p+\epsilon}} \qquad \text{ where, } \eta_{i,i+1}\in [t^n_i,t^n_{i+1}]\\
& \leq \lim_{n\to\infty} \left((\max_i\phi^\prime(\eta_{i,i+1})|\pi^n|)^{\frac{2}{p}-\frac{2}{p+\epsilon}} \times\sum_{\pi^n\cap[0,t]} \left( \bigtriangleup x(t^n_i)\right)^2 \left(\phi^\prime(\eta_{i,i+1}) \bigtriangleup  t_i^n\right)^{\frac{p-2}{p}}\right) \\
& \leq \limsup_{n\to\infty} \left((\max_i\phi^\prime(\eta_{i,i+1})|\pi^n|)^{\frac{2}{p}-\frac{2}{p+\epsilon}}\right) \underbrace{\lim_{n\to\infty}\left( \sum_{\pi^n\cap[0,t]} \left( \bigtriangleup x(t^n_i)\right)^2 \left( \Delta \phi(t^n_{i})\right)^{\frac{p-2}{p}}\right)}_{\langle x \rangle^{(p)}_{\pi}(t)} = 0.
\end{align*}
Since $\phi\in C^1_b([0,T],\R)$, the second equality follows as a consequence of the mean value theorem. 
\par For the other sided equality, consider the case $q=p-\epsilon$, i.e. for all $q<p$, for all $t\in [0,T]$: 
\begin{align*}
\lim_{n\to\infty} 
\mu^{p,q}_{n,\pi}([0,t])& = \lim_{n\to\infty} \sum_{[t^n_i, t^n_{i + 1}] \in \pi^n\cap [0,t]} \left( \bigtriangleup x(t^n_i)\right)^2 \left( \Delta \phi(t^n_i)\right)^{\frac{p-\epsilon-2}{p-\epsilon}} \\
& \geq \lim_{n\to\infty} \left((\max_i \phi'(\eta_{i,i+1})|\pi^n|)^{\left(\frac{2}{p}-\frac{2}{p-\epsilon}\right)} \sum_{\pi^n\cap[0,t]} \left( \bigtriangleup x(t^n_i)\right)^2 \left( \phi'(\eta_{i,i+1})\bigtriangleup t_i^n\right)^{\frac{p-2}{p}}\right) \\
& \geq \left(\liminf_{n\to\infty} \left(\max_i\;\phi'(\eta_{i,i+1})|\pi^n|\right)^{\left(\frac{2}{p}-\frac{2}{p-\epsilon}\right)}\right) \underbrace{\left(\lim_{n\to\infty} \sum_{\pi^n\cap[0,t]} \left( \bigtriangleup x(t^n_i)\right)^2 \left( \Delta\phi(t^n_i)\right)^{\frac{p-2}{p}}\right)}_{\langle x \rangle ^{(p)}_\pi(t)} = \infty,
\end{align*}
the last equality follows from the fact $|\pi^n|\downarrow 0 $, the quantity $\left(\frac{2}{p}-\frac{2}{p-\epsilon}\right)<0$ and the final limit exists and finite as a consequence of $x\in \mathbb S^p_\pi([0,T],\R)$. In second line, there exists $\eta_{i,i+1}\in [t^n_i,t^n_{i+1}]$ using the mean value theorem. This completes part (1). 
\item[(2a)] Since $q>p$ there exists $\epsilon>0$ such that $q=p+\epsilon$. So for all $t \in [0, T]$,
\begin{align*}
0 & \leq \limsup_n\, 
\mu^{p,p+\epsilon}_{n,\pi}([0,t])\\
& \leq \left(\lim_{n\to\infty} (\max_i\phi'(\eta_{i,i+1})|\pi^n|)^{\frac{2}{p}-\frac{2}{p+\epsilon}} \right) \left(\limsup_n \sum_{\pi^n\cap[0,t]} \left( \bigtriangleup x(t^n_i)\right)^2 \left( \Delta\phi( t_{i}^n)\right)^{\frac{p-2}{p}} \right) \\
& \leq C^{\frac{2}{p}-\frac{2}{p+\epsilon}} \left(\lim_{n\to\infty} |\pi^n|^{\frac{2}{p}-\frac{2}{p+\epsilon}} \right) \left(\limsup_n\,   \langle x\rangle_{\pi^n}^{(p)} (t)\right) = 0.
\end{align*}
In the second line, $\phi(t^n_{i+1}) - \phi(t^n_{i}) = \phi'(\eta_{i,i+1}) (t^n_{i+1}-t^n_i)$ for some $\eta_{i,i+1}\in [t^n_i,t^n_{i+1}]$, as a consequence of the mean value theorem. This completes the proof.
\item[(2b)] Since $q<p$ there exists $\epsilon>0$ such that $q=p-\epsilon$. Now; 
\begin{align*}
\limsup_n
\mu^{p,q}_{n,\pi}([0,t])= &\limsup_{n} \sum_{\pi^n\cap [0,t]} \left(|\Delta\phi(t^n_i)|^{1-\frac{2}{p-\epsilon}}\right)(\bigtriangleup x(t^n_i))^2  \\\geq & \left(\lim_{n\to\infty} (C|\pi^n|)^{\frac{2}{p}-\frac{2}{p-\epsilon}} \right) \left(\limsup_n \sum_{ \pi^n\cap [0,t]} \left( \bigtriangleup x(t^n_i)\right)^2 \left( \Delta\phi(t^n_i)\right)^{\frac{p-2}{p}} \right) = \infty.
\end{align*}
In a similar line of argument, (3a) and (3b) follow.
\end{enumerate}
\end{proof}

The following lemma shows that the scaled quadratic variation is invariant under addition by functions which are smoother than the original function.
\begin{lemma}
Let $\pi$ be any sequence of partitions with vanishing mesh and $x\in \mathbb S^p_\pi([0,T],\R)$ for some $p\geq 2$. Then for any $A\in C^0([0,T],\R)$ with $[A]_\pi^{(p)}\equiv 0$ one has:
\[\langle x+A\rangle^{(p)}_\pi \equiv \langle x\rangle^{(p)}_\pi.\]
\end{lemma}
\begin{proof}

Applying Minkowski inequality twice, one gets for all $ t\in[0,T]$,
\begin{align*}
& \left( [x]_{\pi^n}^{(p)}(t)\right)^{1/p} - \left( [A]_{\pi^n}^{(p)}(t)\right)^{1/p} \\
& = \left(\sum_{\pi^n\cap[0,t]} \left|(x+A)(t^n_{j+1}) - (x+A)(t^n_{j}) -A(t^n_{j+1}) +A(t^n_{j}) \right|^p\right)^{1/p} - \left(\sum_{\pi^n\cap[0,t]} \left|A(t^n_{j+1}) - A(t^n_{j}) \right|^p\right)^{1/p} \\
& \leq \left(\sum_{\pi^n\cap[0,t]} \left|(x+A)(t^n_{j+1}) - (x+A)(t^n_{j}) \right|^p\right)^{1/p} = \left( [x+A]_{\pi^n}^{(p)}(t)\right)^{1/p} 
 \leq \left( [x]_{\pi^n}^{(p)}(t)\right)^{1/p} + \left( [A]_{\pi^n}^{(p)}(t)\right)^{1/p}.
\end{align*}
Note that, $\lim_{n\to\infty}[A]_{\pi^n}^{(p)}(t) = [A]_{\pi}^{(p)}(t)\equiv0$ and $\lim_{n\to\infty}[x]_{\pi^n}^{(p)}(t)$ exists and finite, since $x\in \mathbb S^p_\pi([0,T],\R)\subseteq V^p_\pi([0,T],\R)$. Therefore, we conclude that $x+A\in V^p_\pi([0,T],\R)$ and $$[x+A]_{\pi}^{(p)}(t)= [x]_{\pi}^{(p)}(t).$$ 
For $p=2$ the above calculation already concludes the lemma. So without loss of generality, we now assume $p>2$ and drop the indicator $\1$. The scaled quadratic variation of $x+A$ can be written as

\begin{align}\label{eqn1lmp}
\langle x+A\rangle^{(p)}_\pi (t) & = \lim_{n\to\infty} \sum_{\pi^n\cap[0,t]}\left( \Delta [x+A]_{\pi}^{(p)}(t_i^n)\right)^\gamma \left| \bigtriangleup x(t_i^n) + \bigtriangleup A(t_i^n)\right|^2\  \nonumber\\
& = \lim_{n\to\infty} \sum_{\pi^n \cap [0,t]}\left( \Delta [x+A]_{\pi}^{(p)}(t_i^n)\right)^\gamma \left\{ |\bigtriangleup x(t_i^n)|^2 + |\bigtriangleup A(t_i^n)|^2 + 2 \bigtriangleup x(t_i^n)\, \bigtriangleup A(t_i^n)\right\}. 
\end{align}
Note that,
\begin{align*}
& \Delta [x+A]_{\pi}^{(p)}(t_i^n)
= [x+A]_\pi^{(p)}(t_{i+1}^n) - [x+A]_\pi^{(p)}(t_{i}^n)
= [x]_\pi^{(p)}(t_{i+1}^n) - [x]_\pi^{(p)}(t_{i}^n) 
=\bigtriangleup [x]_{\pi}^{(p)}(t_i^n).
\end{align*}
Then,
\begin{align*} 
& \lim_{n\to\infty} \sum_{\pi^n \cap [0,t]}\left( \Delta [x+A]_{\pi}^{(p)}(t_i^n)\right)^\gamma |\bigtriangleup x(t_i^n)|^2 
= \lim_{n\to\infty} \sum_{\pi^n \cap [0,t]}\left( \Delta [x]_{\pi}^{(p)}(t_i^n)\right)^\gamma |\bigtriangleup x(t_i^n)|^2 
= \langle x\rangle_\pi^p(t).
\end{align*}
Again using H\"older's inequality,
\begin{align*}
0\leq & \limsup_{n\to\infty} \sum_{\pi^n \cap [0,t]}\left( \Delta [x+A]_{\pi}^{(p)}(t_i^n)\right)^\gamma |\bigtriangleup A(t_i^n)|^2 = \limsup_{n\to\infty}  \sum_{\pi^n \cap [0,t]}\left( \bigtriangleup [x]_{\pi}^{(p)}(t_i^n) \right)^\gamma |\bigtriangleup A(t_i^n)|^2 \\
& \leq\limsup_{n\to\infty}  \left( \sum_{\pi^n\cap [0,t]}\left( \Delta [x]_{\pi}^{(p)}(t_i^n) \right)^{\gamma\times\frac{p}{p-2}}\right)^{\frac{p-2}{p}} \left( \sum_{\pi^n\cap [0,t]} |\bigtriangleup A(t_i^n)|^{2\times\frac{p}{2}}\right)^{\frac{2}{p}}\\
& = \limsup_{n\to\infty}  \left( \sum_{\pi^n\cap[0,t]}\left( \Delta [x]_{\pi}^{(p)}(t_i^n) \right) \right)^{\gamma} \left( \sum_{\pi^n\cap[0,t]} |\bigtriangleup A(t_i^n)|^{p}\right)^{\frac{2}{p}}\ \text{[as $\gamma:=\frac{p-2}{p}$]} \\
& = \left(\limsup_{n\to\infty} \sum_{\pi^n\cap[0,t]}\left( [x]_\pi^{(p)}(t_{i+1}^n) - [x]_\pi^{(p)}(t_{i}^n)\right) \right)^{\gamma} \left( [A]_\pi^{(p)}(t)\right)^{\frac{2}{p}} \\
& = \left( [x]_\pi^{(p)}(t) - [x]_\pi^{(p)}(0) \right)^{\gamma} \left( [A]_\pi^{(p)}(t)\right)^{\frac{2}{p}}  = \left( [x]_\pi^{(p)}(t)\right)^{\gamma} \left( [A]_\pi^{(p)}(t)\right)^{\frac{2}{p}}  = 0.
\end{align*}
Using Cauchy–Schwarz inequality and the above estimate, we get:
\begin{align*}
0\leq & \left|2 \sum_{\pi^n\cap [0,t]}\left( \Delta [x+A]_{\pi}^{(p)}(t_i^n)\right)^\gamma  \bigtriangleup x(t_i^n)\, \bigtriangleup A(t_i^n)\right| \leq 2 \sum_{\pi^n\cap [0,t]}\left( \Delta [x+A]_{\pi}^{(p)}(t_i^n)\right)^\gamma  |\bigtriangleup x(t_i^n)|\, |\bigtriangleup A(t_i^n)|\\
& \leq 2 \left(\sum_{\pi^n\cap [0,t]}\left( \Delta [x]_{\pi}^{(p)}(t_i^n)\right)^\gamma  |\bigtriangleup x(t_i^n)|^2\right)^{1/2} \times \left(\sum_{\pi^n\cap [0,t]}\left( \Delta [x]_{\pi}^{(p)}(t_i^n)\right)^\gamma  |\bigtriangleup A(t_i^n)|^2\right)^{1/2}\ \\
& \leq 2 \left(\langle x\rangle_{\pi^n}^p\right)^{1/2}\left(\left( [x]_\pi^{(p)}(t)\right)^{\gamma} \left( [A]_{\pi^n}^{(p)}(t)\right)^{2/p}\right)^{1/2}.
\end{align*}
Now taking $\limsup_{n\to \infty}$ we get,
\begin{align*}
0\leq & \limsup_{n\to \infty} \left|2 \sum_{\pi^n\cap [0,t]}\left( \Delta [x+A]_{\pi}^{(p)}(t_i^n)\right)^\gamma  \bigtriangleup x(t_i^n)\, \bigtriangleup A(t_i^n)\right| \\& \hspace{3cm} \leq  2 \left(\langle x\rangle_{\pi}^p\right)^{1/2}\left(\left( [x]_\pi^{(p)}(t)\right)^{\gamma} \left( [A]_{\pi}^{(p)}(t)\right)^{2/p}\right)^{1/2} = 0
\end{align*}
By \eqref{eqn1lmp} we now conclude the proof.
\end{proof}

\section{Scaled quadratic variation under smooth transformation}\label{sec:3}
One of the crucial properties of $p$-th variation is that for any partition sequence with vanishing mesh, $V^p_\pi$ invariant under smooth transformations. i.e. if $x\in V^p_\pi([0,T],\R)$ then for any `smooth' $f$, $f\circ x \in V^p_\pi([0,T],\R)$. In fact when $f\in C^2$ and $x\in V^p_\pi([0,T],\R)$, one can write the $p$-th variation of $f\circ x$ as a Riemann–Stieltjes integral involving $p$-th variation of $x$ (see \cite[Theorem 2.1]{perkowski2019}, \cite{ananova2017}).

In this section, we provide such a transformation formula for scaled quadratic variation.

\begin{theorem}[Transformation rule for pathwise scaled quadratic
variation]\label{Thm.ito.isomerty}
Let $p\geq2$, $f \in C^{2}(\R,\R)$ and $\pi$ be a partition sequence of $[0,T]$ with vanishing mesh. Let $x\in \mathbb S^p_\pi([0,T],\R)\cap C^\alpha([0,T],\R) $ for some $\alpha> \frac{\sqrt{1+\frac{4}{p}}-1}{2}$, i.e. 
$ \lim_n \langle x\rangle^{(p)}_{\pi^n}$ exists and finite.  Then
\[\langle f \circ x\rangle^{(p)}_\pi(t) = \int_0^t |f^\prime\circ x(u)|^p\, d \langle x\rangle^{(p)}_\pi(u).\]
\end{theorem}
\begin{proof}
For notational convenience, we write
$\phi(t):=[x]^{(p)}_\pi(t)$, $\phi^f(t):=[f\circ x]^{(p)}_\pi(t)$ and $\gamma = \frac{p-2}{p}$. Since we work with $p \geq 2$, we drop the $\1$.
From the definition of scaled quadratic variation with index $p$,
\begin{align}
& \langle f \circ x\rangle^{(p)}_\pi(t) 
 = \lim_{n\to\infty} \sum_{\pi^n\cap [0,t]} \left( \phi^f(t_{i+1}^n) - \phi^f(t_{i}^n)\right)^\gamma \left| (f\circ x)(t_{i+1}^n) - (f\circ x)(t_{i}^n)\right|^2   \nonumber\\
& \overset{\text{\cite[Theorem 3.3]{ruhong2024}}}{=} \lim_{n\to\infty}  \sum_{\pi^n\cap [0,t]}  \left( \int_0^{t^n_{i+1}} \left\{|f^\prime(x(u))|^p - \int_0^{t^n_{i}} |f^\prime(x(u))|^p\right\} d\phi(u)\right)^\gamma \left| (f\circ x)(t_{i+1}^n) - (f\circ x)(t_{i}^n)\right|^2 
\nonumber\\
& = \lim_{n\to\infty}  \sum_{\pi^n\cap [0,t]}  \left(  \int_{t_i^n}^{t^n_{i+1}} |f^\prime(x(u))|^p\,   d\phi(u)\right)^\gamma \left| (f\circ x)(t_{i+1}^n) - (f\circ x)(t_{i}^n)\right|^2.\  \label{Eq: Mvt}
\end{align}
As a consequence of the mean value theorem \cite[ Problem 1.1.26]{KaczorNowak3}, for all $i \in \{0,1,2,\cdots, N(\pi^n)-1\}$ there exists $v_{i}\in [t_i^n,t_{i+1}^n]$ such that $\int_{t_i^n}^{t^n_{i+1}} |f^\prime(x(u))|^p\,   d\phi(u) = |f^\prime(x(v_i))|^p\, \left(\Delta \phi(t_{i}^n)\right)$. So the above formula can be simplified further as follows. 
\begin{align}
& = \lim_{n\to\infty}  \sum_{\pi^n\cap [0,t]} \left( |f^\prime(x(v_i))|^p\, \left( \Delta \phi(t_{i}^n) \right)\right)^\gamma\left| (f\circ x)(t_{i+1}^n) - (f\circ x)(t_{i}^n)\right|^2\  \nonumber\\
& = \lim_{n\to\infty}  \sum_{\pi^n\cap [0,t]}  |f^\prime(x(v_i))|^{p\gamma}\, 
 \left(\Delta \phi(t_{i}^n) \right)^\gamma \left| (f\circ x)(t_{i+1}^n) - (f\circ x)(t_{i}^n)\right|^2. \label{bom1}
\end{align}
We will now show that \eqref{bom1} is equal to  $\int_0^t |f^\prime(x(u))|^p\, d\langle x\rangle_\pi^{(p)}(u)$ completing the proof. As a consequence of Taylor expansion there exists  $\eta_i\in[x(t_{i}^n)\wedge x(t_{i+1}^n), x(t_{i}^n)\vee x(t_{i+1}^n)]$ such that 
\begin{align}\label{bom2}
& \left| (f\circ x)(t_{i+1}^n) - (f\circ x)(t_{i}^n)\right|^2 \nonumber\\
& = \left| f^\prime (x(t_i^n))(\bigtriangleup x (t^n_i)) + \frac{1}{2} f^{\prime\prime} (\eta_i)(\bigtriangleup x (t^n_i))^2 \right|^2, \nonumber\\
& = \left(f^\prime (x(t_i^n))\right)^2(\bigtriangleup x (t^n_i))^2 + \frac{1}{4} \left(f^{\prime\prime} (\eta_i)\right)^2(\bigtriangleup x (t^n_i))^4 + f^\prime (x(t_i^n)) f^{\prime\prime} (\eta_i) (\bigtriangleup x (t^n_i))^3.
\end{align}
Now, from \eqref{bom1} and \eqref{bom2} consider the term
\begin{align*}
& \lim_{n\to\infty}  \sum_{\pi^n\cap [0,t]}  |f^\prime(x(v_i))|^{p\gamma}\, \left(\Delta \phi(t_{i}^n) \right)^\gamma \left(f^\prime (x(t_i^n))\right)^2(\bigtriangleup x (t^n_i))^2 \\
& = \lim_{n\to\infty} \sum_{\pi^n\cap [0,t]}  |f^\prime(x(v_i))|^{p-2} \left(f^\prime (x(t_i^n))\right)^2\, \left(\Delta[x]_\pi^{(p)} (t_i^n) \right)^\gamma (\bigtriangleup x (t^n_i))^2 \\
& = \int_0^t |f^\prime(x(u))|^p\, d\langle x\rangle_\pi^{(p)}(u).
\end{align*}
Again, from \eqref{bom1} and \eqref{bom2} consider the term
\begin{align*}
& \left| \limsup_{n\to\infty}  \sum_{\pi^n\cap [0,t]}  |f^\prime(x(v_i))|^{p\gamma}\, 
 \left(\Delta \phi(t_{i}^n) \right)^\gamma \left|f^\prime (x(t_i^n)) f^{\prime\prime} (\eta_i) (\bigtriangleup x (t^n_i))^3\right| \right| \\
 & \leq \left\|f^\prime \circ x\right\|_\infty^{p-1}  \left\|f^{\prime\prime}\circ x\right\|_\infty \left\{ \lim_{n\to\infty}  \sum_{\pi^n\cap [0,t]} \left|\Delta[x]_\pi^{(p)} (t_i^n) \right|^\gamma (\bigtriangleup x (t^n_i))^2 \right\} \\
 & \qquad\qquad\hspace{7 cm} \times\left\{ \lim_{n\to\infty} \sup_{[t_i^n,t^n_{i+1}]\in\pi^n\cap [0,t]} |\bigtriangleup x (t^n_i)|\right\}
  \\
 & = 0.
\end{align*}
In the above argument, we use 1) $  \sup_{[u,v]\in\pi^n\cap [0,t]} |x(v)-x(u)|\to0$, as a result of $|\pi^n|\downarrow0$, uniform continuity of $x$, and 2) finiteness of remaining terms as a consequence of the assumptions on $x$ and $f$.
\par Similarly, from \eqref{bom1} and \eqref{bom2} consider the term
\begin{align*}
& \left| \limsup_{n\to\infty}  \sum_{\pi^n\cap [0,t]}  |f^\prime(x(v_i))|^{p\gamma}\, 
 \left(\Delta \phi(t_{i}^n) \right)^\gamma \frac{1}{4} \left(f^{\prime\prime} (\eta_i)\right)^2(\bigtriangleup x (t^n_i))^4\right| \\
 & \leq \frac{1}{4}\,  \left\|f^\prime\circ x\right\|_\infty^{p-2}  \left\|f^{\prime\prime}\circ x\right\|_\infty^2 \left\{ \lim_{n\to\infty}  \sum_{\pi^n\cap [0,t]} \left|\Delta[x]_\pi^{(p)} (t_i^n) \right|^\gamma (\bigtriangleup x (t^n_i))^2 \right\} \\
 & \quad\hspace{7cm}\quad \times \left\{ \lim_{n\to\infty} \sup_{t_i^n\in\pi^n\cap [0,t]} \left|\bigtriangleup x (t^n_i)\right|^2\right\} \\
 & = 0,
\end{align*}
where the argument is similar to the above. This concludes the proof. 
\end{proof}

\begin{remark}[The case $1<p<2$ of Theorem \ref{Thm.ito.isomerty}]\label{rem: p<2}
For $p<2$ the exponent $\gamma=\frac{p-2}{p}$ is negative, and the proof of 
Theorem~\ref{Thm.ito.isomerty} uses $p\ge2$ in exactly two places: the bound 
$|f'(x(v^n_i))|^{p-2}\le\|f'\circ x\|_\infty^{p-2}$ in the two remainder 
estimates, and the uniform continuity of $y\mapsto|y|^{p-2}$ on 
$f'(x([0,T]))$ in the convergence of the main term. Both fail if $f'\circ x$ 
vanishes somewhere, since $|y|^{p-2}\to\infty$ as $y\to0$. To extend Theorem \ref{Thm.ito.isomerty} for the case of $1<p<2$ we need  additional non-degeneracy 
condition
\[
\delta:=\inf_{u\in[0,T]}|f'(x(u))|>0 .
\]
Indeed under the extra assumption $|f'(x(v^n_i))|^{p-2}\le\delta^{p-2}$, the two remainder terms  are bounded by $\delta^{p-2}\|f'\circ x\|_\infty\|f''\circ x\|_\infty\, \langle x\rangle^{(p)}_{\pi^n}(t)\,\mathbf{osc}(x,\pi_n)$ and 
$\tfrac14\delta^{p-2}\|f''\circ x\|_\infty^2\,\langle x\rangle^{(p)}_{\pi^n}(t)\, \mathbf{osc}(x,\pi_n)^2$ respectively. 
Note also that in this case  $\Delta\phi^f(t^n_i) = \int_{t^n_i}^{t^n_{i+1}}|f'(x(u))|^p\,d\phi(u) \ge\delta^p\,(\Delta\phi(t^n_i))$, so the indicators 
$\mathbbm 1_{\Delta\phi^f(t^n_i)>0}$ and $\mathbbm 1_{\Delta\phi(t^n_i)>0}$ 
coincide and the convention \eqref{eq: convention} can be applied to $\Delta\phi(t^n_i)$ throughout the proof. 
\end{remark}
\begin{remark}
One advantage of the above transformation rule (in the functional calculus setup) is that the path in consideration can be `arbitrarily rough' (i.e. all $p\in (1,\infty)$). In particular, as $p\to \infty$ the H\"older condition $\alpha > \frac{1}{2}\left(\left(1 + \frac{4}{p} \right)^{\frac{1}{2}} - 1\right)$ in Theorem \ref{Thm.ito.isomerty} becomes $C^{0+}$. Alternative approaches to these problems, such as the controlled Rough Path approach, have an inherent restriction of $\alpha \in \left(\frac{1}{3}, 1\right]$ (See, for example, \cite[Theorem 3.5]{leahy-EJP25}).
\end{remark}

\begin{remark}
Under some additional mild conditions, in the spirit of \cite{ruhong2024}, the above theorem may be generalised to a functional set-up. In the absence of a mean value theorem-type result for an integral setup with a functional integrand, one can use a first-order Taylor expansion in \eqref{Eq: Mvt}. 
\end{remark}

\section{Examples}\label{sec:4}
For a one-dimensional fractional Brownian motion with Hurst exponent $H$, as a consequence of results in \cite{gladyshev1961} and \cite{pratelli2011}, the scaled quadratic variation with index $\frac1H$, exists and is finite. We nonetheless explicitly calculate the scaled quadratic variation for a two-dimensional fractional Brownian path.
\begin{example}[Scaled quadratic variation for Correlated fractional Brownian Motion] \label{ex:fbm_covariation} 
Let $H\in (0,1), \rho \in [-1,1]$, $\gamma = 1-2H$ and consider two independent fractional Brownian motions $x^1,\tilde x$ with Hurst exponent $H$ on the same probability space $(\Omega,\mathcal F, \mathbb P)$. Now for all $\omega\in\Omega$ construct the process $x^2(\omega) := \rho x^1(\omega) + \sqrt{1-\rho^2}\tilde{x}(\omega)$. So $x^2$ is simply another fractional Brownian motion on the same probability space $(\Omega,\mathcal F, \mathbb P)$ with correlation coefficient between $x^1$ and $x^2$ is $\rho$. Now we will calculate the scaled quadratic variation of the two-dimensional continuous process $x = (x^1, x^2) \in C([0,T], \mathbb{R}^2)$, which is a two-dimensional fractional Brownian motion with Hurst parameter $H$ and cross-correlation $\rho \in [-1,1]$, so that $\mathbb{E}[x^1(t)x^2(t)] = \rho t^{2H}$ for all $t \in [0,T]$. 

Consider the sequence of dyadic partitions $\pi = (\pi_n)_{n \ge 1}$ of $[0,T]$ defined by $t_j^n = \frac{j T}{2^n}$ for $j = 0, 1, \dots, 2^n$, with mesh size $|\pi_n| = \frac{T}{2^n} \to 0$. 
Along the dyadic sequence of partitions $(\pi_n)_n$ (note that this is subsequence of uniform partition), the $1/H$-th variation of each component path $x^k$ ($k=1,2$) is known to converge deterministically (see Mishura \cite[Section 1.18]{mishura2008book}, Pratelli \cite{pratelli2011}) to:
\[
[x^k]_\pi^{(1/H)}(t) = c_H t, \quad \text{where } c_H := \mathbb{E}\left[|Z|^{1/H}\right] \text{ with } Z \sim \mathcal{N}(0,1).
\]
Note here that the assumptions of Lemma \ref{SQV-recovers-pth-var} are not satisfied, and scaled quadratic variation of index $\frac{1}{H}$ for $x^k$ turns out to be a constant multiple of $\frac{1}{H}$-th variation of the same path.
\par Now, following the notation of Definition \ref{def:scaled_covariation} the variation increments are constant across intervals: $\Delta [x^k]_\pi^{(1/H)}(t_j^n) = c_H \Delta t_j^n = c_H \frac{T}{2^n} > 0$. Substituting these increments into the scaled quadratic covariation measure $\mu_{n,\pi}^{(1,2)}$, we obtain for $t \in [0,T]$:
\begin{align*}
\mu_{n,\pi}^{(1,2)}([0,t]) &= \sum_{j=0}^{\lfloor 2^n t / T \rfloor - 1} \left( \Delta [x^1]_\pi^{(1/H)}(t_j^n) \Delta [x^2]_\pi^{(1/H)}(t_j^n) \right)^{\gamma/2} \bigtriangleup x^1 (t^n_j) \, \bigtriangleup x^2 (t^n_j) \\
&= \sum_{j=0}^{\lfloor 2^n t / T \rfloor - 1} \left( c_H \frac{T}{2^n} \right)^{1-2H} \bigtriangleup x^1 (t^n_j) \, \bigtriangleup x^2 (t^n_j).
\end{align*}
Using the decomposition $x^2 = \rho x^1 + \sqrt{1-\rho^2}\tilde{x}$, and the fact $x^1$ and $\tilde{x}$ are two independent fBM with Hurst index $H$ on the same probability space, we get
\begin{align*}
\mu_{n,\pi}^{(1,2)}([0,t]) = &\rho \, \underbrace{\left(c_H\frac{T}{2^n}\right)^{1-2H} 
               \sum_{j=0}^{\lfloor 2^nt/T \rfloor - 1}(\bigtriangleup x^1 (t^n_j))^2 }_{A_n(t)} 
    \\ &+ \sqrt{1-\rho^2} \, \underbrace{\left(c_H\frac{T}{2^n}\right)^{1-2H} 
    \sum_{j=0}^{\lfloor 2^nt/T \rfloor - 1} (\bigtriangleup x^1 (t^n_j))(\bigtriangleup \tilde{x} (t^n_j))}_{M_n(t)}.    
\end{align*}
By the Gladyshev theorem \cite[Theorem~1]{gladyshev1961} (see also 
\cite[Section~1.18]{mishura2008book}) we have
\[
    A_n(t) \xrightarrow{n\to\infty} c_H^{1-2H} \cdot t \quad \text{almost surely } .
\]
For the term $M_n(t)$, we will use $4ab = (a+b)^2 - (a-b)^2$ formula. 
\begin{align*}
    4M_n(t) = \left(c_H\frac{T}{2^n}\right)^{1-2H} \sum_{j=0}^{\lfloor 2^nt/T \rfloor-1} 
    \left[
    \left((x^1 +\tilde{x})(t^n_{j+1}) - (x^1 +\tilde{x})(t^n_{j})\right)^2 
    - 
    \left((x^1 -\tilde{x})(t^n_{j+1}) - (x^1 -\tilde{x})(t^n_{j})\right)^2
    \right].
\end{align*}
Since $x^1$ and $\tilde{x}$ are independent fBMs with Hurst index $H$, the processes $x^1 + \tilde{x}$ and $x^1 - \tilde{x}$ are both Gaussian  processes on the original probability space with covariance function $(|t|^{2H}+|s|^{2H}-|t-s|^{2H})$. So we can now apply the Gladyshev theorem \cite[Theorem~1]{gladyshev1961} individually to the processes $x^1 + \tilde{x}$ and $x^1 - \tilde{x}$ resulting in
\begin{align*}
    M_n(t) \xrightarrow[]{n\to \infty} 0 \quad \forall t\in [0,T] \qquad \text{almost surely.}
\end{align*}
Hence, the matrix-valued scaled quadratic variation process of the two-dimensional fBM is given by;
\[
\langle x \rangle_\pi^{(1/H)}(t) = c_H^{1-2H} t \begin{pmatrix} 1 & \rho \\ \rho & 1 \end{pmatrix}, \quad \text{ for all } t \in [0,T] \text{ almost surely} .
\]
\end{example}

Below we give examples of functions (beyond the Gaussian setup, hence not covered by Gladyshev's theorem \cite[Theorem 1]{gladyshev1961}) for which the (pathwise) scaled quadratic variation exists and is finite.

\begin{example}\label{Ex:Takagi}
Consider the generalised Takagi class of functions (for details see \cite[Equation (2)]{Misura2019} and the references therein) as follows, for all $H\in (0,1)$: 
\begin{align*}
\mathcal X^H= \left\{x\in C^0([0,1],\R) \;\big|\; x(t) =\sum_{m=0}^\infty 2^{m(\frac{1}{2}-H)}\sum^{2^m-1}_{k=0}
\theta_{m,k}e_{m,k}(t) \text{ for coefficients } \theta_{m,k}\in \{-1, +1\}\right\}.
\end{align*}
Then the scaled quadratic variation with index $\gamma = 1-2H$ along the sequence $\mathbb T$ of dyadic partitions $\mathbb T_n$ of $[0,1]$ exists and is finite.  To see this fix $H\in (0,1)$. Take any $x\in \mathcal X^H$, then
\begin{align*}
\langle x \rangle^{(1/H)}_{\mathbb T}(t) &= \lim_{n\to \infty} \sum_{\mathbb T^n\cap[0,t]} ([x]^{(p)}_\mathbb T (t^n_{i+1}) - [x]^{(p)}_\mathbb T (t^n_{i}))^{1-2H} (\bigtriangleup x (t^n_i))^2\, \1_{\Delta [x]^{(p)}_\mathbb T (t^n_{i})}\\
&= \lim_{n\to \infty} \sum_{\mathbb T^n\cap[0,t]} ( C_H t^n_{i+1}-C_H t^n_{i})^{1-2H} (\bigtriangleup x (t^n_i))^2,
\end{align*} 
where the last equality follows from \cite[Theorem 2.1]{Misura2019}.  Please see \cite{Misura2019} for further details on the constant $ C_H$. Since $\mathbb T$ is dyadic, $t^n_{i+1} - t^n_{i} = \frac{1}{2^n}$. So the above can be reduced further to (for notational convenience, assume $t=1$):
\begin{align*}
\langle x \rangle^{(1/H)}_{\mathbb T}(1) & = \lim_{n\to\infty} C_H^\gamma \frac{1}{2^{n(1-2H)}} [x]_{\mathbb T^n} (1) = \lim_{n\to\infty} C_H^\gamma \frac{1}{2^{n(1-2H)}} \frac{1}{2^n}\sum_{m=0}^{n-1}\sum_{k=0}^{2^m-1}2^{m(1-2H)}\theta_{m,k}^2\\
&= \lim_{n\to\infty} C_H^\gamma \frac{1}{2^{n(1-2H)}} \frac{1}{2^n}\sum_{m=0}^{n-1} 2^{2m(1-H)} = \lim_{n\to\infty} C_H^\gamma \frac{1}{2^{n(1-2H)}}\frac{1}{2^n} \frac{2^{2n(1-H)}-1}{2^{2-2H}-1} = \frac{ C_H^{1-2H}}{2^{2-2H}-1}.
\end{align*}
Following a similar line of argument in \cite{das2021} we have $x\in \mathbb S^p_{\mathbb T}([0,1],\R)$ and $\langle x \rangle^{(\gamma)}_{\mathbb T}(t) =  \frac{  C_H^{1-2H}}{2^{2-2H}-1}t$ completing the proof. 
\end{example}
\begin{example}\label{Ex:fakeFBM}
Consider the fake fractional Brownian motion class of functions (for details see \cite[Sections 6 \& 7]{bayraktar2025}) as follows, for all $H\in (0,1)$:
\begin{align*}
Y^H (t) =\sum^\infty_{m=0}\sum_{k} \theta^{Y, H}_{m,k}e_{m,k}(t),
\end{align*}
where the Schauder coefficients $\theta^{Y, H}_{m,k}$ satisfy the assumptions  of \cite[Corollary 7.4]{bayraktar2025}.   
Since the $1/H$-th variation of this fake process is linear, as shown in \cite[Corollary 7.4]{bayraktar2025}, the proof follows in a similar line of argument to Example \ref{Ex:Takagi}. So the scaled quadratic variation along the partition sequence $\pi$ with index $\gamma = 1-2H$ exists and is finite, where $1/H$ is an even integer.
\end{example}
In line with the above two examples, we note that the existence and finiteness of scaled quadratic variation does not imply that this is the same as $p$-th variation.

In the spirit of normalized $p$-th variation \cite[Section~4.3]{das2022a} one can introduce a notion of normalized (pathwise) scaled quadratic variation, which can be used as an estimator of variation index (see \cite[Section~4.4]{das2022a} and the references therein). We leave this as future work. 

\textbf{Acknowledgement}: We thank the anonymous reviewer for several constructive feedback, including the suggestion to focus on the self-normalization part of the results, leading to a significant improvement in the presentation. We would like to thank Fang Rui Lim for the counterexample in Example \ref{counter.example.1}, Rama Cont for constructive feedback and Torstein Nilssen for insightful discussion. This research is partially supported by the Jointly Funded Bilateral Mobility Program - Partnership Collaboration Awards between Indian Institute of Technology Madras and King’s College London. Suprio Bhar acknowledges the support of the Matrics grant MTR/2021/000517 from the Science and Engineering Research Board (Department of Science \& Technology, Government of India).

\bibliographystyle{siam}
\bibliography{pathwise1}

@article{Misura2019,
title = {{On (signed) Takagi–Landsberg functions: $p$th variation, maximum, and modulus of continuity}},
journal = {Journal of Mathematical Analysis and Applications},
volume = {473},
number = {1},
pages = {258-272},
year = {2019},
issn = {0022-247X},
doi = {https://doi.org/10.1016/j.jmaa.2018.12.047},
url = {https://www.sciencedirect.com/science/article/pii/S0022247X18310862},
author = {Yuliya Mishura and Alexander Schied}}

@article{das2021,
 title = {Quadratic variation along refining partitions: Constructions and examples},
journal = {Journal of Mathematical Analysis and Applications},
volume = {512},
number = {2},
pages = {126 -- 173},
year = {2022},
issn = {0022-247X},
author = {Rama Cont and Purba Das}
}

@article{ananova2017,
	Author = {Anna Ananova and Rama Cont},
	Doi = {http://dx.doi.org/10.1016/j.matpur.2016.10.004},
	Issn = {0021-7824},
	Journal = {Journal de Math{\'e}matiques Pures et Appliqu{\'e}es},
	Pages = {737--757},
	Title = {Pathwise integration with respect to paths of finite quadratic variation},
	Url = {http://www.sciencedirect.com/science/article/pii/S0021782416301155},
	Year = {2017},
	volume = {107},
	number={6}
}

@incollection{pratelli2011,
author="Pratelli, Maurizio",
title="A Remark on the $1/{H}$-Variation of the Fractional {Brownian} Motion",
bookTitle="S{\'e}minaire de Probabilit{\'e}s XLIII",
year="2011",
publisher="Springer Berlin Heidelberg",
address="Berlin, Heidelberg",
pages="215--219",
isbn="978-3-642-15217-7",
doi="10.1007/978-3-642-15217-7_8",
url="https://doi.org/10.1007/978-3-642-15217-7_8"
}

@article{perkowski2019,
	Author = {Cont, Rama and Perkowski, Nicolas},
	Fjournal = {Transactions of the American Mathematical Society},
	Issn = {},
	Journal = {Transactions of the American Mathematical Society},
	Mrclass = {},
	Mrnumber = {},
	Mrreviewer = {161-161},
	Pages = {161--186},
	Title = {Pathwise integration and change of variable formulas for continuous paths with arbitrary regularity},
	Volume = {6},
	Year = {2019}}

@incollection{follmer1981,
  title={Calcul d'{I}t{\^o} sans probabilit{\'e}s},
  author={F{\"o}llmer, Hans},
  booktitle={S{\'e}minaire de Probabilit{\'e}s XV 1979/80: Avec table g{\'e}n{\'e}rale des expos{\'e}s de 1966/67 {\`a} 1978/79},
  pages={143--150},
  year={2006},
  publisher={Springer}
}

@article{schied2016,
	Author = {Schied, Alexander},
	journal = {Journal of Mathematical Analysis and Applications},
	Issn = {0022-247X},
	Mrclass = {26A27 (26A15 60G17 60H05)},
	Mrnumber = {3398747},
	Number = {2},
	Pages = {974--990},
	Title = {On a class of generalized {T}akagi functions with linear pathwise quadratic variation},
	Volume = {433},
	Year = {2016}
}

@book{tudor2013,
	Author = {Tudor, Ciprian A},
	Isbn = {},
	Mrclass = {Book},
	Mrnumber = {},
	Note = {Part of the Probability and Its Applications book series (PIA)},
	Pages = {271},
	Publisher = {Springer},
	Title = {Analysis of Variations for Self-similar Processes: A Stochastic Calculus Approach},
	Url = {https://link.springer.com/book/10.1007%2F978-3-319-00936-0},
	Year = {2013}}

@article{ruhong2024,
  title={{Fractional It\^o calculus}},
  author={Cont, Rama and Jin, Ruhong},
  journal={Transactions of the American Mathematical Society, Series B},
  volume={11},
  number={22},
  pages={727--761},
  year={2024}
}

@article{davis2018,
author = "Davis, Mark and Ob{\l}{\'o}j, Jan and Siorpaes, Pietro",
doi = "10.1214/16-AIHP792",
fjournal = "Annales de l'Institut Henri Poincaré, Probabilités et Statistiques",
journal = "Ann. Inst. H. Poincaré Probab. Statist.",
month = "02",
number = "1",
pages = "1--21",
publisher = "Institut Henri Poincaré",
title = "Pathwise stochastic calculus with local times",
url = "https://doi.org/10.1214/16-AIHP792",
volume = "54",
year = "2018"
}

@book{freedman,
	Author = {Freedman, D.},
	Title = {{Brownian} Motion and Diffusion },
	Year = {1983},
	publisher = {Springer}}

@article{gladyshev1961,
  title={A New Limit Theorem for Stochastic Processes with {G}aussian Increments},
  author={Gladyshev, EG},
  journal={Theory of Probability \& Its Applications},
  volume={6},
  number={1},
  pages={52--61},
  year={1961},
  publisher={SIAM}
}

@book {KaczorNowak3,
    AUTHOR = {Kaczor, W. J. and Nowak, M. T.},
     TITLE = {Problems in mathematical analysis. {III}},
    SERIES = {Student Mathematical Library},
    VOLUME = {21},
      NOTE = {Integration},
 PUBLISHER = {American Mathematical Society, Providence, RI},
      YEAR = {2003},
     PAGES = {x+356},
      ISBN = {0-8218-3298-0},
   MRCLASS = {28-01 (00A07 26-01)},
  MRNUMBER = {1998016},
MRREVIEWER = {B.\ Crstici},
       DOI = {10.1090/stml/021},
       URL = {https://doi.org/10.1090/stml/021},
}

@article{bayraktar2025,
  title={H{\"o}lder regularity and roughness: construction and examples},
  author={Bayraktar, Erhan and Das, Purba and Kim, Donghan},
  journal={Bernoulli},
  volume={31},
  number={2},
  pages={1084--1113},
  year={2025},
  publisher={Bernoulli Society for Mathematical Statistics and Probability}
}

@article {leahy-EJP25,
    AUTHOR = {Leahy, James-Michael and Nilssen, Torstein},
     TITLE = {Scaled quadratic variation for controlled rough paths and
              parameter estimation of fractional diffusions},
   JOURNAL = {Electron. J. Probab. (updated version available in arXiv preprint \url{https://arxiv.org/abs/2401.09299v4})},
  FJOURNAL = {Electronic Journal of Probability},
    VOLUME = {30},
      YEAR = {2025},
     PAGES = {},
      ISSN = {1083-6489},
   MRCLASS = {60H10 (60F15 60G17 60G22 60L20 62F12 62M09)},
  MRNUMBER = {4888157},
       DOI = {10.1214/25-ejp1308},
       URL = {https://doi.org/10.1214/25-ejp1308},
}

@article {MR2965922,
    AUTHOR = {Kubilius, K. and Mishura, Y.},
     TITLE = {The rate of convergence of {H}urst index estimate for the
              stochastic differential equation},
   JOURNAL = {Stochastic Process. Appl.},
  FJOURNAL = {Stochastic Processes and their Applications},
    VOLUME = {122},
      YEAR = {2012},
    NUMBER = {11},
     PAGES = {3718--3739},
      ISSN = {0304-4149,1879-209X},
   MRCLASS = {60G22 (60H10 62F12)},
  MRNUMBER = {2965922},
MRREVIEWER = {B.\ L. S. Prakasa Rao},
       DOI = {10.1016/j.spa.2012.06.011},
       URL = {https://doi.org/10.1016/j.spa.2012.06.011},
}

@phdthesis{das2022a,
  edition = {},
  number = {},
  journal = {},
  pages = {},
  publisher = {University of Oxford},
  school = {University of Oxford},
  title = {Roughness properties of paths and signals},
  volume = {},
  author = {Das, P},
  editor = {},
  year = {2022},
  series = {}
}

@article{daskim2025,
title = {On isomorphism of the space of continuous functions with finite $p$-th variation along a partition sequence},
journal = {Journal de Mathématiques Pures et Appliquées},
volume = {203},
note         = {Paper No.\ 103753},
year = {2025},
issn = {0021-7824},
doi = {https://doi.org/10.1016/j.matpur.2025.103753},
url = {https://www.sciencedirect.com/science/article/pii/S0021782425000972},
author = {Purba Das and Donghan Kim}
}

@book {jacod-protter-discretization,
    AUTHOR = {Jacod, Jean and Protter, Philip},
     TITLE = {Discretization of processes},
    SERIES = {Stochastic Modelling and Applied Probability},
    VOLUME = {67},
 PUBLISHER = {Springer, Heidelberg},
      YEAR = {2012},
     PAGES = {xiv+596},
      ISBN = {978-3-642-24126-0},
   MRCLASS = {60F05 (60G44 60G51 60G57 60H10 60H35 60J25 60J75)},
  MRNUMBER = {2859096},
MRREVIEWER = {Dominique\ L\'epingle},
       DOI = {10.1007/978-3-642-24127-7},
       URL = {https://doi.org/10.1007/978-3-642-24127-7},
}

@book {mishura2008book,
    AUTHOR = {Mishura, Yuliya S.},
     TITLE = {Stochastic calculus for fractional {B}rownian motion and
              related processes},
    SERIES = {Lecture Notes in Mathematics},
    VOLUME = {1929},
 PUBLISHER = {Springer-Verlag, Berlin},
      YEAR = {2008},
     PAGES = {xviii+393},
      ISBN = {978-3-540-75872-3},
   MRCLASS = {60G15 (60-02 60H05 60H10 62M09 62M20 62P05 91B28)},
  MRNUMBER = {2378138},
MRREVIEWER = {Ivan\ Nourdin},
       DOI = {10.1007/978-3-540-75873-0},
       URL = {https://doi.org/10.1007/978-3-540-75873-0},
}

@article{Baxter1956,
 ISSN = {00029939, 10886826},
 URL = {http://www.jstor.org/stable/2032765},
 author = {Glen Baxter},
 journal = {Proceedings of the American Mathematical Society},
 number = {3},
 pages = {522--527},
 publisher = {American Mathematical Society},
 title = {A Strong Limit Theorem for {G}aussian Processes},
 urldate = {2026-09-10},
 volume = {7},
 year = {1956}
}

@book{mishura2017-book,
  title={Parameter estimation in fractional diffusion models},
  author={Kubilius, K{\k{e}}stutis and Mishura, Yuliya and Ralchenko, Kostiantyn},
  volume={8},
  year={2017},
  publisher={Springer}
}

@article{kubilius2010,
  title={On the convergence rates of {G}ladyshev’s {H}urst index estimator},
  author={Kubilius, Kestutis and Melichov, Dmitrij},
  journal={Nonlinear Analysis: Modelling and Control},
  volume={15},
  number={4},
  pages={445--450},
  year={2010}
}

@article{han-Schied2025,
  title={The roughness exponent and its model-free estimation},
  author={Han, Xiyue and Schied, Alexander},
  journal={The Annals of Applied Probability},
  volume={35},
  number={2},
  pages={1049--1082},
  year={2025},
  publisher={Institute of Mathematical Statistics}
}

\end{document}